\documentclass{birkjour}

 \newtheorem{thm}{Theorem}[section]
 
 \newtheorem{lem}[thm]{Lemma}
 
 \theoremstyle{definition}
 
 \theoremstyle{remark}
 \newtheorem{rem}[thm]{Remark}
 
 \numberwithin{equation}{section}

\begin{document}

\title[Oberbeck-Boussinesq equations]
 {On the  regularity of weak solution for the Oberbeck-Boussinesq equations}
\author[B. Climent-Ezquerra]{Blanca Climent-Ezquerra}
\address{Departamento de Ecuaciones Diferenciales y An\'alisis Num\'erico\br
Universidad de Sevilla\br
Sevilla, Espa\~na}
\email{bcliment@us.es}
\author[E. Ortega-Torres]{Elva Ortega-Torres}
\address{Departamento de Matem\'aticas\br
Universidad Cat\'olica del Norte\br
Antofagasta, Chile}
\email{eortega@ucn.cl}
\author[M. Rojas-Medar]{Marko Rojas-Medar}
\address{Departamento de Matem\'atica\br
Universidad de Tarapac\'a\br
Arica, Chile}
\email{mmedar@academicos.uta.cl}
\date{}

\begin{abstract}
We prove new regularity criteria of the Prodi-Serrin type with weak Lebesgue integrability
in both space and time for a viscous active chemical fluid in a bounded domain.
\end{abstract}

\maketitle

\section{Introduction}

In this paper, we show some regularity results for the equations describing the motion of a viscous
chemical fluid in a bounded domain, $\Omega \subset \mathbb{R}^3$, with smooth boundary, $\partial \Omega$,
over a time interval, $[0,T)$, $0<T\leq \infty$. Specifically, we consider the following approximations
of Oberbeck-Boussinesq type (see \cite{joseph}):
\begin{equation}\label{In.uno}
\left \{ \begin{array}{l}
{\boldsymbol u}_t+ ({\boldsymbol u}\cdot\nabla){\boldsymbol u} -\mu \Delta{\boldsymbol u} +\nabla p
=\alpha (\theta+\varphi){\boldsymbol g}+{\boldsymbol f},\\
\theta_t + ({\boldsymbol u}\cdot \nabla) \theta -\kappa_1 \Delta \theta = \ell,\\
\varphi_t + ({\boldsymbol u}\cdot\nabla)\varphi -\kappa_2\Delta\varphi = h,\\
\mbox{div}\,{\boldsymbol u} =0,
\end{array}\right.
\end{equation}
together with the following boundary and initial conditions:
\begin{equation}\label{In.cin}
{\boldsymbol u}(x,t)=0, \ \theta(x,t)= 0, \ \varphi(x,t)=0 \mbox{ on } \partial \Omega \times (0,T).
\end{equation}
\begin{equation}\label{In.cua}
{\boldsymbol u}(x,0)= {\boldsymbol u}_0(x), \ \theta(x,0)=\theta_0(x), \ \varphi(x,0)= \varphi_0  \mbox{ in }\Omega.
\end{equation}
The unknowns are the functions ${\boldsymbol u}(x,t) \in \mathbb{R}^3$, $\theta(x,t) \in \mathbb{R}$, $\varphi (x,t)\in
\mathbb{R}$ and $p(x,t) \in \mathbb{R}$ denoting the velocity vector, temperature, concentration of
material in the liquid and pressure at time $t \in [0,T)$, at point $x\in \Omega$. Moreover,
${\boldsymbol f}(x,t), {\boldsymbol g} (x,t) \in \mathbb{R}^3$, $\ell(x,t), h(x,t) \in \mathbb{R}$ are known external sources,
$\mu>0$ is the viscosity of fluid, $\kappa_1$ and $\kappa_2$ are the thermal and solute diffusion,
respectively.  ${\boldsymbol u}_0, \theta_0$  and $ \varphi_0$ are functions given on the variable $x \in \Omega$.
The non-homogeneous case for $\theta$ and $\varphi$ in the initial conditions can be treated using an
appropriate lifting and only the obvious changes in the statement of the results are required.

The expressions $\Delta$, $\nabla$ and {\rm div} denote the gradient, Laplacian and divergence operators, respectively;
the $i^{th}$ component of $({\boldsymbol u}\cdot\nabla){\boldsymbol u}$  is given by $[({\boldsymbol u}
\cdot \nabla){\boldsymbol u}]_i= \sum_{j=1}^3 u_j \frac{\partial u_i}{\partial x_j}$
and $({\boldsymbol u}\cdot\nabla)\phi = \sum_{j=1}^3 u_j \frac{\partial \phi}{\partial x_j}$ for $\phi = \theta$ or $\varphi$.
For the derivation and physical discussion of equation (\ref{In.uno})-(\ref{In.cua}) see, e.g., Joseph \cite{joseph}.

We note that this model of fluid includes as a particular case the classical Navier-Stokes, which has
been extensively studied (see, for instance, the classic books by Ladyzhenskaya \cite{Lady}, Lions \cite{Lions}
and Temam \cite{Temam}). It also includes the classical Boussinesq problem (or Bernard's problem) in the case where
chemical reactions are absent: $\varphi \equiv 0$. The existence of various solution notions for this system has been
studied, see \cite{shinbroth}, \cite{korenev}, \cite{canon}, \cite{morimoto1}, \cite{hishida}, \cite{morimoto2} and
\cite{komo} and the references therein.

The study of the uniqueness of weak solution for three-dimensional Navier-Stokes equations is more difficult, see for example \cite{olga}.
Uniqueness is proved fairly easily for sufficiently regular solutions while in the class of weak solutions, the problem of uniqueness
(or non-uniqueness) remains open. A pioneering work in this direction was given by Serrin \cite{Serrin}, who proved that a weak solution
$({\boldsymbol u}, p)$ is regular if the velocity ${\boldsymbol u}$ satisfies some suitable additional conditions (see next section). 
Afterwards, many authors have extended Serrin's result, see e.g.  \cite{kaniel}, \cite{von}, \cite{Veiga2}, \cite{Bjorland}, 
\cite{Bosia}, \cite{Pata}, \cite{Sohr}, etc.

The model considered in this work was studied by Belov and Kapitonov, who proved the stability of the solutions in \cite{belov}. Rojas-Medar
and Lorca established the existence, regularity and uniqueness of the solutions by means of the Galerkin spectral method in
\cite{lorca1}, \cite{lorca2}. In \cite{morimoto2}, Morimoto obtained results analogous to those in Serrin \cite{Serrin} involving
conditions on ${\boldsymbol u}$ and $\theta$ for the Boussinesq approximation. This results can be extended to the model considered in 
this paper, using the results of weak solutions obtained in \cite{Sebastian}. Note that this extension involves 
conditions of the Prodi-Serrin type for ${\boldsymbol u}$, $\theta$ and $\varphi$.  Here, we obtain regularity results 
assuming only suitable conditions on ${\boldsymbol u}$. This work is inspired by the article \cite{Bosia}
(see also \cite{miguel}) .

The paper is organized as follows:  In section 2, we show the notation and some preliminary results. The two theorems presented in this
paper are stated in section 3 and in the following sections their proofs are given.

\section{Preliminaries}

\subsection{Notation}
Let $\Omega$ be an open subset of $\mathbb{R}^n$. Boldface letters will be used for
vectorial spaces, for instance ${\textbf{L}^2(\Omega)}=L^2(\Omega)^n$. Sometimes, the notation will be
abridged, if $X=X(\Omega)$ is a space of
functions defined in the open set $\Omega$, we denote by $L^p(X)$ the
Banach space $L^p(0,T;X)$. The $L^p$-norm is denoted as $\|\cdot\|_p$ (in particular
$\|\cdot\|_2=\|\cdot\|$). The inner product
of $L^2(\Omega)$ is denoted by $(\cdot,\cdot)$.

As usual, we define $\mathcal{V}=\{{\boldsymbol v} \in {\bf C}_0^\infty(\Omega)\, : \, \mbox{div}\, 
{\boldsymbol v}=0\}$ and the spaces
\begin{align*}
{\bf H}&= \mbox{ the closure of }\mathcal{V} \mbox{ in } {\bf L}^2(\Omega),\\
{\bf V}&= \mbox{ the closure of }\mathcal{V} \mbox{ in } {\bf H}^1(\Omega).
\end{align*}
We denote by $P$ the orthogonal projection onto ${\bf H}$ and $A=-P \Delta$ the usual Stokes operator with
domain $D(A)={\bf H}^2(\Omega)\cap {\bf V}$.

From now on, $C>0$ will denote different constants, depending only on the  data of the problem. When
necessary, we emphasize that the constants may have different values using the notation $C_1$, $C_2$,
and so on.

\subsection{Weak Lebesgue or Marcinkiewicz spaces}
 Given a measurable function $f$ on $\Omega$, the distribution function of $f$ is the function $d_f$, defined on $[0,\infty)$ by
$d_f(\alpha)=m(\{x \in \Omega; |f(x)|>\alpha\})$. For $1\leq p<\infty$, the weak-$L^p(\Omega)$ space is denoted by $L^{p,\infty}(\Omega)$
and is defined as the set of measurable functions $f$ such that
\begin{align*}
\| f\|_{p,\infty} & = \sup_{\alpha}\left \{\alpha \cdot d_f(\alpha)^{1/p}; \alpha>0\right\}
\end{align*}
is finite. It is possible to check that $\| \lambda f\|_{p,\infty}=|\lambda| \| f\|_{p, \infty}$, $\| f\|_{p, \infty}=0$
if and only if $f=0$, and $\| f+g\|_{p,\infty} \leq 2(\| f\|_{p,\infty}+ \| g\|_{p,\infty})$.  Hence, $L^{p,\infty}(\Omega)$
is a quasinormed linear space (for more details see Grafakos \cite{Graf}, Kufner et al. \cite{kufner}).

The weak-$L^\infty(\Omega)$ space is by definition the space $L^\infty(\Omega)= L^{\infty,\infty}
(\Omega).$ The space $L^{p,\infty}(\Omega)$ is larger than
the usual $L^p(\Omega)$ space. This follows directly from Chebyshev's inequality
\begin{align*}
\alpha^p d_f(\alpha) & \leq \int_{\{x; |f(x)|>\alpha\}} |f(x)|^pdx.
\end{align*}
More precisely, the inclusion $L^p(\Omega) \subset L^{p,\infty}(\Omega)$ is strict. For example, if $x_0 \in \Omega$,
we have $|x-x_0|^{-n/p} \in L^{p,\infty}(\Omega)$, but $|x-x_0|\notin L^{p}(\Omega).$

 We will also use the following lemma (Lemma 2.18.2 in \cite{kufner}):
\begin{lem}\label{Lem.un} If $f\in L^p(\Omega)$, then
\begin{align*}
\int_\Omega\| f(x)\|^p\,dx &= p \int_0^\infty \sigma^{p-1}d_f(\sigma)\,d\sigma.
\end{align*}
\end{lem}
\subsection{Previous results}

For the Navier-Stokes equations it is well known that if ${\boldsymbol u}(0)\in {\bf H}$ and ${\boldsymbol f} \in L^2(0,T;
{\bf L}^2(\Omega))$, there exists at least one weak solution 
\[{\boldsymbol u}\in L^\infty(0,T; {\bf H})\cap L^2(0,T; {\bf V}),\]
which is called the {\it Leray-Hopf} solution. If, in addition ${\boldsymbol u}(0) \in {\bf V}$, there exists 
$T_* \in (0, \infty]$ such that the Navier-Stokes equations admit a unique strong solution
\[{\boldsymbol u} \in L^\infty(0,T; {\bf V}) \cap L^2(0,T; {\bf H}^2(\Omega)),\]
provided that $T<T^*$. Therefore, it is natural to ask under what conditions do you have a global in
time strong solution for initial data  ${\boldsymbol u}(0)\in {\bf V}$, without any restriction on the size of its norm.
According to a result of Prodi \cite{Prody} and Serrin \cite{Serrin}, such a global strong solution for
the Navier-Stokes equations exists if
\begin{align*}
&{\boldsymbol u} \in   L^s(0,T; {\bf L}^r(\Omega)),
\end{align*}
where $(r,s)$ is a Prodi-Serrin pair, that is, $r \in (3, \infty], s \in [2, \infty)$ and $\frac{3}{r}+\frac{2}{s}=1$.
Moreover, using weak Lebesgue spaces, it has been showed that the weak solution of the Navier-Stokes
equations is strong on $[0,T]$ if either
\begin{align*}
& {\boldsymbol u}\in L^s(0,T; {\bf L}^{r,\infty}(\Omega)) \ \mbox{ or } \
{\boldsymbol u} \in L^{s,\infty}(0,T; {\bf L}^{r,\infty}(\Omega))
\end{align*}
and it exists a constant $\Gamma=\Gamma(r,s, \Omega)>0$ such that
 $\|{\boldsymbol u}\|_{L^{s,\infty}({\bf L}^{r,\infty})}\leq \Gamma\, \mu^{1-1/s}$ where $\Gamma>0$ is a small constant
 (see \cite{Bjorland}, \cite{Bosia}, \cite{Sohr}).

 The following lemma of Bosia {\it et al} \cite{Bosia} will be necessary in our demonstrations:
\begin{lem}\label{Lem.uno} Let ${\boldsymbol u}\in {\bf V} \cap {\bf H}^2(\Omega)$ and  $(r,s)$ be a Prodi-Serrin 
pair, that is, $r \in (3,\infty], s \in [2, \infty)$ and $\frac{3}{r}+\frac{2}{s}=1$. Then
\begin{align*}
\|({\boldsymbol u}\cdot\nabla){\boldsymbol u} & \|\leq C_r \|{\boldsymbol u}\|_{r, \infty} \|\nabla{\boldsymbol u}\|^{\frac{2}{s}}
\|A {\boldsymbol u}\|^{\frac{s-2}{s}}.
\end{align*}
\end{lem}
\section{The main results}
Our purpose here is to extend the above regularity criteria of the Navier-Stokes equations to problem
(\ref{In.uno}) and provide short proofs of  the results.

We assume that a weak solution $({\boldsymbol u}, \theta, \varphi)$ of (\ref{In.uno}) satisfies
\begin{align*}
&{\boldsymbol u} \in L^{\infty}(0,T; {\bf H})\cap L^2(0,T; {\bf V}),\\
&\theta, \varphi \in L^\infty(0,T; L^2(\Omega))\cap L^2(0,T; H^1_0(\Omega)).
\end{align*}
\begin{thm} Let ${\boldsymbol f}\in L^2(0,T; \L^2(\Omega))$, ${\boldsymbol g} \in L^{\infty}(0,T; 
{\bf L}^3(\Omega))$ and $\ell, h \in L^2(0,T; L^2(\Omega))$. Let $({\boldsymbol u}, \theta, \varphi)$  
be a weak solution of (\ref{In.uno}) with initial data ${\boldsymbol u}(0) \in {\bf V}$ and $ \theta(0), $ 
$\varphi(0) \in H_0^1(\Omega)$. If
\begin{align*}
&{\boldsymbol u} \in L^s(0,T; {\bf L}^{r,\infty}(\Omega))
\end{align*}
for some Prodi-Serrin pair $(r,s)$, then $({\boldsymbol u},\theta, \varphi)$ remains strong on $[0,T]$.
\label{Th.uno}
\end{thm}

\begin{thm}\label{Th.dos} Let ${\boldsymbol f} \in L^2(0,T; {\bf L}^2(\Omega))$, ${\boldsymbol g} 
\in L^{\infty}(0,T; {\bf L}^3(\Omega))$ and $\ell, h \in L^2(0,T; L^2(\Omega))$. Let $({\boldsymbol u}, 
\theta, \varphi)$  be a weak solution of (\ref{In.uno}) with initial data ${\boldsymbol u}(0) \in {\bf V},$ 
$ \theta(0), $ $\varphi(0) \in H_0^1(\Omega)$. If there exists  a constant $\Gamma>0$, depending only 
on $r$ and $\Omega$, such that
\begin{align*}
{\boldsymbol u} \in L^{s,\infty}(0,T; {\bf L}^{r,\infty}(\Omega))
\mbox{ with } \|{\boldsymbol u}\|_{L^{s,\infty}({\bf L}^{r,\infty})}& \leq \Gamma \,\mu^{(s-1)/s}
\end{align*}
for some Prodi-Serrin pair $(r,s)$, then $({\boldsymbol u},\theta, \varphi)$ remains strong on $[0,T]$.
\end{thm}

\begin{rem} The constant $\Gamma$ is  $\Gamma = (2C^\epsilon_2 C_1)^{-1/s}$, where $C_1$ and
$C_2$ are given in the proof of Theorem \ref{Th.dos} and $\epsilon>0$ is appropriately chosen.
\end{rem}

\section{Proof of Theorem \ref{Th.uno}}
By applying the operator $P$ in (\ref{In.uno}) and taking $A{\boldsymbol u}$ as test function in the 
weak formulation of (\ref{In.uno}) we get
\begin{align}\label{P.tre}
\frac{1}{2}\frac{d}{dt}\|\nabla{\boldsymbol u}\|^2+\mu\|A{\boldsymbol u}\|^2&
= -( P[({\boldsymbol u}\cdot\nabla ){\boldsymbol u}], A {\boldsymbol u}) 
+ \alpha( P(\theta +\varphi){\boldsymbol g}, A{\boldsymbol u})\nonumber\\
&\ \ \ + ({\boldsymbol f},A{\boldsymbol u}).
\end{align}
Lemma \ref{Lem.uno} provides the following estimate
\begin{align*}
|( P[({\boldsymbol u}\cdot \nabla ){\boldsymbol u}], A {\boldsymbol u}) |& \leq \|({\boldsymbol u}
\cdot \nabla){\boldsymbol u}\|\|A {\boldsymbol u}\|\\
&\leq C_s \|{\boldsymbol u}\|_{r,\infty}\|\nabla{\boldsymbol u}\|^{\frac{2}{s}}\|A{\boldsymbol u}\|^{\frac{2(s-1)}{s}}\\
&\leq  \frac{C_1}{2}\mu^{1-s} \|{\boldsymbol u}\|_{r,\infty}^s \|\nabla{\boldsymbol u}\|^2
+ \frac{\mu}{6} \|A{\boldsymbol u}\|^2,
\end{align*}
where $ \frac{C_1}{2}=\frac{C_s^s}{s}\left [\frac{6(s-1)}{s}\right]^{s-1}$.
By using the Holder's and Young's inequalities we obtain
\begin{align*}
|\alpha ( P \theta{\boldsymbol g}, A {\boldsymbol u}))|& \leq \alpha \|(\theta +\varphi){\boldsymbol g}\|
\|A {\boldsymbol u}\|\\
&\leq \frac{6}{\mu} \alpha^2\|(\theta+\varphi){\boldsymbol g}\|^2+ \frac{\mu}{6}\|A{\boldsymbol u}\|^2\\
&\leq  C_2\frac{6}{\mu} \alpha^2 \|{\boldsymbol g}\|_{L^\infty({\bf L}^3)}(\|
\nabla\theta\|^2 + \|\nabla\varphi\|^2) +\frac{\mu}{6} \|A {\boldsymbol u}\|^2.
\end{align*}
where $C_2$ is the embedding constant $H^1 \hookrightarrow L^6$ and,
\[|({\boldsymbol f}, A{\boldsymbol u})| \leq \frac{6}{\mu} \|{\boldsymbol f}\|^2 + \frac{\mu}{6} 
\|A {\boldsymbol u}\|^2.\]
From the above estimates and from (\ref{P.tre}), we get the following differential inequality:
\begin{align}\label{P.cua}
\frac{d}{dt}\|\nabla{\boldsymbol u}\|^2 +\mu\|A{\boldsymbol u}\|^2 &\leq  C_1\mu^{1-s} \|{\boldsymbol u}\|_{r,\infty}^s
\|\nabla{\boldsymbol u}\|^2\nonumber\\
& \ \ \ + \tilde{C}(\|\nabla\theta\|^2 +\|\nabla \varphi\|^2 )+\frac{\mu}{3}\|{\boldsymbol f}\|^2,
\end{align}
where $\tilde{C} = 2C_2 \frac{6}{\mu} \alpha^2 \|{\boldsymbol g}\|_{L^\infty({\bf L}^3)}$.

Setting $\phi(t)=\|\nabla{\boldsymbol u}(t)\|^2$ and $k(t)=\|{\boldsymbol u}(t)\|_{r,\infty}^s$, inequality 
(\ref{P.cua}) implies
\begin{align}\label{Pa.dos}
\phi(t)&\leq  \phi(0)+ C_1 \int_0^t k(\sigma)\phi(\sigma)d\sigma
+\tilde{C}\int_0^t (\|\nabla \theta(\sigma)\|^2+\|\nabla \varphi(\sigma)\|^2)d\sigma\nonumber\\
& \ \ \ + \frac{\mu}{3} \int_0^t \|{\boldsymbol f}(\sigma) \|^2 d \sigma\nonumber\\
&\leq  M + C_1\int_0^t k(\sigma) \phi(\sigma)d\sigma,
\end{align}
where $ M= \varphi(0)+ \tilde{C}( \|\theta\|_{L^2(H_0^1)} +\|\varphi\|_{L^2(H_0^1)})
+ \frac{\mu}{3} \|{\boldsymbol f}\|_{L^2({\bf L}^2)}  <\infty$.

By Gronwall's inequality we conclude that
\[\|\nabla{\boldsymbol u}(t) \|^2=\varphi(t)\leq  M \exp \bigg(C_1\int_0^t k(\sigma)d\sigma\bigg).\]
Hence,
\begin{equation}\label{Pa.tre}
{\boldsymbol u} \in L^\infty(0,T; {\bf V}).
\end{equation}
Similarly,  by taking the inner product of the second equation of (\ref{In.uno}) with $-\Delta \theta$ we have that
\begin{equation}\label{P.dos}
\frac{1}{2}\frac{d}{dt}\|\nabla \theta\|^2 + \kappa_1\|\Delta \theta\|^2
= -(\ell,\Delta\theta)+ ( ({\boldsymbol u}\cdot \nabla) \theta, \Delta \theta).
\end{equation}
By using the Holder's and Young's inequalities, we have
\[|-(\ell, \Delta \theta)|\leq \|\ell\|\|\Delta \theta\|\leq C \|\ell\|^2+\frac{\kappa_1}{4}\|\Delta
\theta\|^2.\]
Owing to the Gagliardo-Niremberg, Poincar\'e, Young inequalities and  (\ref{Pa.tre}), we obtain
\begin{align}\label{P.sie}
|(({\boldsymbol u}\cdot\nabla)\theta, \Delta\theta)|&\leq \|{\boldsymbol u}\|_6 \| \nabla \theta \|_3 
\| \Delta \theta\|\nonumber\\
&\leq C\|\nabla{\boldsymbol u}\| (\|\nabla\theta\|^{\frac{1}{2}}\|\Delta\theta\|^{\frac{1}{2}})
\|\Delta\theta\|\nonumber\\
&\leq C \|\nabla \theta\|^2+ \frac{\kappa_1}{4}\|\Delta \theta\|^2.
\end{align}
Using the estimates above in (\ref{P.dos}) we get
\begin{align}\label{P.cin}
|\nabla\theta(t)\|^2 &\leq \|\nabla \theta(0)\|^2 +C \int_0^t
\|\nabla\theta(\sigma)\|^2d\sigma +C\int_0^t \|\ell(\sigma)\|^2d\sigma\nonumber\\
& \leq N +C \int_0^t \|\nabla \theta(\sigma)\|^2d\sigma,
\end{align}
where $N=\|\nabla \theta(0)\|^2 +\|\ell\|_{L^2(0,T; L^2(\Omega) }<\infty$. By the Gronwall's
inequality, we conclude that
\begin{equation}\label{Pa.cua}
\theta \in L^\infty(0,T; H^1_0(\Omega)).
\end{equation}
Analogously,
\begin{equation}\label{Pa.cua1}
\varphi \in L^\infty(0,T; H^1_0(\Omega)).
\end{equation}
can be shown.

We then assume that the solution remains strong only on $[0,T')$, with $T'<T$. It follows from
(\ref{Pa.tre}), (\ref{Pa.cua}) and (\ref{Pa.cua1}) that $\|\nabla{\boldsymbol u}\|^2+ \|\nabla\theta\|^2
+ \|\nabla\varphi\|^2$ remains bounded on $[0,T]$.  But if $[0,T')$ is the maximal interval
for the existence of the strong solution, then $\|\nabla{\boldsymbol u}\|+\|\nabla\theta\|+\|\nabla\varphi\|
 \to \infty$ as $t \to T'^{-}$ and we reach a contradiction. Hence, the solution must be strong
 on $[0,T]$.

\section{Proof of Theorem \ref{Th.dos}}
We use an argument similar to the one used in Bosia {\it et al} \cite{Bosia} and Pata \cite{Pata}: Given
$ 0<\epsilon<1$, let $s_\epsilon=s+\epsilon(4-s)$ and  $r_\epsilon$ such that $(r_\epsilon, s_\epsilon)$
is a Prodi-Serrin pair. By using the following  interpolation inequality
\[\|{\boldsymbol u}\|_{r_\epsilon,\infty}^{s_\epsilon}\leq \|{\boldsymbol u}\|_{r,\infty}^{s(1-\epsilon)}
\|{\boldsymbol u}\|_{6,\infty}^{4\epsilon}
 \leq C_2^\varepsilon \|{\boldsymbol u}\|_{r, \infty}^{s(1-\epsilon)}\|\nabla{\boldsymbol u}\|^{4\epsilon},\]
in (\ref{P.cua}) for $s_\epsilon$ instead of $s$, it follows that
\[\frac{d}{dt} \|\nabla{\boldsymbol u}\|^2\leq C_3\mu^{1-s} \|{\boldsymbol u}\|_{r,\infty}^{s(1-\epsilon)}
\|\nabla{\boldsymbol u}\|^{2+4\epsilon}+ C (\|\nabla\theta\|^2+\|\nabla\varphi\|^2+ \|{\boldsymbol f}\|^2).\]
In this inequalities, $C_2^\epsilon$ depends only on $\Omega$ and $\epsilon$,
$C_3=C_2^\epsilon C_1$ and $C_1$ depends on $(r_\epsilon, s_\epsilon)$.
Hence, setting $\phi(t)=\| \nabla{\boldsymbol u}(t)\|^2$ and integrating in $[0,t]$ we obtain
\begin{equation}\label{Pa.cin}
\phi(t)\leq R + C_3\mu^{1-s} \int_0^t\|{\boldsymbol u}(\sigma)\|_{r, \infty}^{s(1-\epsilon)}
\phi(\sigma)^{1+2\epsilon} d\sigma,
\end{equation}
where $R= \phi(0)+ C (\|\theta\|_{L^2( H_0^1)}+\|\varphi\|_{L^2(H_0^1)}+\|{\boldsymbol f}\|_{L^2({\bf L}^2)})
< \infty.$

We denote $k(t)= \|{\boldsymbol u}\|_{r,\infty}^{s}$ and
\[\psi(t)= R+ C_3\mu^{1-s} \int_0^t k^{1-\epsilon}(\sigma) \phi^{1+2\epsilon}(\sigma)d\sigma
\ \mbox{ for } \ t \in [0,T].\]
Observe that $k \in L^{1, \infty}(0,T)$ and $\phi\leq \psi$. From (\ref{Pa.cin}), we have that
\begin{align*}
\psi'(t)&= C_3\mu^{1-s} k^{1-\epsilon}(t) \phi^{1+2 \epsilon}(t)\\
&\leq C_3 \mu^{1-s} k^{1-\epsilon}(t) \psi^{1+2 \epsilon}(t).
\end{align*}
Integrating the last inequality in $[0,t]$, we obtain
\begin{equation}\label{M.bos}
\psi(0)^{-2\epsilon}- \psi(t)^{-2\epsilon}
\leq 2\, C_3\, \mu^{1-s}\epsilon\int_0^t k^{1-\epsilon}(\sigma)d\sigma.
\end{equation}
On the other hand, owing to Lemma \ref{Lem.un},
\begin{align}\label{M.uno}
\epsilon\int_0^T k^{1-\epsilon}(\sigma)d\sigma &=\epsilon(1-\epsilon)\int_0^\infty
\sigma^{-\epsilon} m\{\tau \in [0,T]; k(\tau)> \sigma\} d\sigma\nonumber\\
&\leq \epsilon T +\epsilon (1-\epsilon)\| k \|_{1, \infty}\int_1^\infty
\sigma^{-(1+\epsilon)}d\sigma\nonumber\\
&= \epsilon T+(1-\epsilon)\| k\|_{1, \infty}.
\end{align}
By using (\ref{M.uno}) in (\ref{M.bos}),
\begin{equation}\label{M.dos}
\psi(0)^{-2\epsilon}- \psi(t)^{-2\epsilon} \leq 2C_3\mu^{1-s} \epsilon T
+ 2 C_3\mu^{1-s} (1-\epsilon)\| k\|_{1,\infty}.
\end{equation}
Let $\delta \in (0,1/3)$ such that $2\,C_3\,\mu^{1-s}\|k\|_{1, \infty}<1-3\delta$. Choosing
$\epsilon>0$ sufficiently small  so that $1-\delta<\psi^{-2\epsilon}(0)$ and $2 C_3 \mu^{1-s}
\epsilon T <\delta$, we have $\psi(t)\leq \delta^{-1/(2\epsilon)}$ for all $t \in [0,T].$
Therefore, if $ \| k \|_{1, \infty}<1/(2C_3\mu^{1-s}  )$, then $\phi=\| \nabla{\boldsymbol u}\|^2$ is
bounded on $[0,T]$.
Arguing as in the proof of Theorem \ref{Th.uno}, we can conclude that $\|\nabla\theta\|$ and
$\|\nabla \varphi\|$ are bounded on $[0,T]$ (see (\ref{Pa.cua})).


\subsection*{Acknowledgment}
M. Rojas-Medar is partially supported by MATH-AMSUD project 21-MATH-03 (CTMicrAAPDEs), CAPES-PRINT 
88887.311962/2018-00 Brazil, Project UTA-Mayor, 753-20, Universidad de Tarapac\'a-Chile.

\end{document}